\documentclass[12pt, letterpaper]{article}
\usepackage{amsfonts}
\usepackage{amsmath}

\newcommand{\sect}[1]{\addtocounter{section}{1}~\newline
{\large \bf \thesection.  #1}\\}

\renewcommand{\Re}{{\rm Re}}
\renewcommand{\Im}{{\rm Im}}

\newcommand{\rank}{{\rm rank}}

\newcommand{\infpoint}{\{\infty\}}

\newcommand{\cR}{\mathcal{R}}

\newcommand{\cX}{\mathcal{X}}
\newcommand{\cY}{\mathcal{Y}}

\newcommand{\cD}{\mathcal{D}}
\newcommand{\cK}{\mathcal{K}}

\newcommand{\cC}{\mathcal{C}}

\newcommand{\cF}{\mathcal{F}}

\newcommand{\bi}{{\mathbf{i}}}

\newcommand{\bR}{\mathbb{R}}
\newcommand{\bC}{\mathbb{C}}

\newcommand{\bHM}{\mathbb{H}\mathbb{M}}

\newcommand{\pv}{p}
\newcommand{\pd}{\cD}

\newcommand{\be}{\begin{equation}}
\newcommand{\ee}[1]{\label{#1} \end{equation}}

\newcommand{\eq}[1]{(\ref{#1})}

\newcommand{\vect}[2]{
\left( \begin{array}{c} #1\\ #2 \end{array} \right) }
\newcommand{\matr}[4]
{\left( \begin{array}{cc} #1 & #2\\ #3 & #4
\end{array} \right) }

\newtheorem{theorem}{Theorem}
\newtheorem{corollary}{Corollary}

\newcommand{\comment}[1]{}

\title{ Parameter-Dependent S-Procedure And Yakubovich
Lemma\\[1ex]
{Sergei V. Gusev}\\[1ex]
{\small St.~Petersburg State University, Russia\\[1.5ex]
 {\it E-mail:} gusev@ieee.org\\
{\it Web page:} http://www.math.spbu.ru/user/gusev}
}

\begin{document}
\begin{center}
{\large \bf Parameter-Dependent S-Procedure And Yakubovich
Lemma}\\
\end{center}
\begin{center}
{\large Sergei V. Gusev}\\
\end{center}
{\normalsize\it Department of Mathematics and Mechanics,}
{\normalsize\it St Petersburg State University } \newline
{\normalsize\it 2 Bibliotechnaya sq., Peterhof,
St.Petersburg, 198904, Russia}\newline
E-mail address: gusev@ieee.org\newline
Web Page: http://www.math.spbu.ru/user/gusev\\[3ex]

\begin{center} { \bf Abstract} \end{center}
The paper considers  a linear matrix
inequality (LMI) that depends on a parameter
varying in a compact topological space.
It turns out that if a strict LMI continuously depends
on a parameter and is feasible for any value of that parameter,
then it has a solution which continuously depends on the parameter.
The result holds true for LMIs that arise
in S-procedure and Yakubovich lemma.
It is shown that the  LMI  which is polynomially
dependent on a vector of parameters  can be reduced to
a parameter-independent LMI of a higher dimension.
The result is based on the recent generalization of
Yakubovich lemma proposed by Iwasaki and Hara and
another generalization formulated in this paper.
The problem of positivity verification for a
non-SOS polynomial of two variables
is considered as an example.
To illustrate control applications,
a method  of parameter-dependent
Lyapunov function construction is proposed  for
nonlinear systems with parametric uncertainty.
\\[1em]

\noindent {\sl Comment:}
The paper was presented at the 5th Russian-Swedish Control Conference,
Lund, Sweden, 29--30 August 2006.\\


\thispagestyle{empty}
\newpage
\begin{flushright}
{\it  To Vladimir A. Yakubovich\\
on the occasion of his eightieth birthday}
\end{flushright}
\sect{Introduction}

The S-procedure losslessness theorem and
the Kalman-Yakubovich-Popov (KYP) lemma
are important mathematical tools of modern control theory.
Both statements deal with inequalities
of certain quadratic forms.
We consider a situation when these quadratic forms depend
continuously on a parameter that varies in a compact topological
space.

First, we study solutions of a parameter-dependent
linear matrix inequality (LMI).
It turns out that if a strict LMI continuously depends
on a parameter and is feasible, then it has  a solution
that continuously depends on the parameter.
From this general statement it follows that
if S-procedure for strict inequalities
is lossless for any value of the parameter,
then the  Lagrange multipliers
can be chosen as continuous functions of the parameter.

We also consider parameter-dependent generalizations of one
statement of KYP lemma, namely the equivalence of strict
frequency-domain inequality and the strict linear matrix inequality.
(This result was first proven
by V.A.~Yakubovich~\cite{Yakubovich.1962}
and was called Yakubovich lemma by
R.~Kalman~\cite{Kalman.1963} and
S.~Lefschetz~\cite{Lefschetz.1965}.)

It is proven that
if the matrices in the frequency-domain inequality
continuously depend on the parameter, then
there exists a solution of the LMI,
which continuously depends on the parameter.
The result also holds true  for the generalized KYP
lemma proposed by T.~Iwasaki and S.~Hara~\cite{Iwasaki.Hara.AC.2005}
as well as for a new version of KYP lemma with matrix
frequency-domain inequality.
The different result concerning a parameter-dependent
version of KYP lemma was obtained by
A.~L.~Likhtarnikov~\cite{Likhtarnikov.2001}.
In~\cite{Likhtarnikov.2001}  solutions of Lur'e equation
are  considered that
are analytic functions of a parameter.

Using parameter-dependent Yakubovich lemma
we show that if an LMI  polynomially depends on the parameters,
then it can be transformed into an LMI of a larger dimension
that does not depend on these parameters.
As an illustration of this result,
we consider minimization of a polynomial
of several variables in a bounded domain.
The proposed method can  be used to verify the positivity
of polynomials that cannot be represented
as a sum of squares.

The results are also applied  to the construction
of a parameter-dependent Lyapunov function for
stability analysis of nonlinear systems
with parametric uncertainty.

\sect{Parameter-dependent LMI}

Consider the parameter-dependent LMI
in $h \in \bR^{l}$
\be
L(\pv, h) > 0
\ee{L>0}
where $p \in \cD$ is a parameter,
$L : \pd\times \bR^{l} \to \bHM_n$  is affine
with respect to $h.$
Hereafter $\bHM_n$ denotes the space
of Hermitian matrices of dimension $n.$
\begin{theorem} Suppose that
$\cD$ is a compact topological space,
the map $L$ is continuous and is affine
with respect to $h,$
and that  LMI \eq{L>0} is feasible for any $p \in \cD;$ then
there exists a continuous function $h(.) : \pd \to \bR^{l}$
that satisfies inequality
\be
L(\pv, h(\pv)) > 0 \; \; \forall \pv \in \pd.
\ee{L>0forallp}

The set of all continuous  $h(.)$ satisfying
\eq{L>0forallp} is open in $\cC(\pd, \bR^l).$
\end{theorem}
Hereafter $\cC(\cX, \cY)$ denotes the normed space of
continuous functions from the topological space $\cX$
to the normed vector space $\cY.$

\begin{corollary}
If $\cF$ is a dense subset of $\cC(\pd, \bR^{l}),$
then there exists  $h(.) \in \cF$  satisfying \eq{L>0forallp}.
\end{corollary}

\begin{corollary} \label{cor2}
If $\pd \subset \bR^{k},$
then there exists a polynomial $h(\pv)$ satisfying
\eq{L>0forallp}.
\end{corollary}

Suppose that $h(p)$ is a polynomial of
a vector variable $p \in \bR^{k}.$
Denote  the vector
of coefficients of $h(\pv)$
by $\hat{h} \in \bR^{\hat{k}}.$
Let $\hat{L} : \pd\times\bR^{\hat{k}} \to \bHM_n$
be a map that is defined
by substitution of the polynomial $h(p)$ into $L.$
From corollary~\ref{cor2} it follows that
\be
\hat{L}(\pv, \hat{h}) >0 \; \; \forall \pv \in \pd.
\ee{hatL>0}

This way the construction of
the parameter-dependent solution $h(\pv)$
of \eq{L>0} is reduced to a search of
a constant vector $\hat{h}$ that
satisfies the system of inequalities \eq{hatL>0}
indexed by $\pv.$
The most interesting case is the one when this system includes
infinite number of LMIs.
In section~7 we show that under certain assumptions
on $\cD$ this system can be reduced to
a single finite-dimensional LMI.


\sect{Parameter-dependent S-procedure}

Consider two maps
$F_0 : \cD\times\cX \to \bR$
and $F : \cD\times\cX \to \bR^l,$
where $\cD$ is a compact space,
$\cX$ is either $\bC^n$ or $\bR^n.$
Let $\cK $ be closed convex proper cone in $\bR^l.$
The dual cone is denoted by
$\cK^+ = \{ z \in \bR^l | \langle z, y \rangle \ge 0 \; \forall y \in
\cK\}.$ 
Consider two statements:\\

I. $\forall \pv \in \pd \;\;\; F(\pv, x)
\in \cK \Rightarrow F_0(\pv, x) > 0.$\\

II. $\exists h : \pd \to \cK^+$ such that
\be
F_0(\pv, x) > \langle h(\pv), F(\pv, x) \rangle \;
\forall \pv \in \pd, x \in \cX.
\ee{F0>hF}

It is clear that II implies I.
In the case when parameter $p$ is absent,
the operation of  replacement of the condition I by
the condition II is known as S-procedure.
In the presence of parameter we call this method
a parameter-dependent S-procedure.
The S-procedure is said to be lossless
if I is equivalent to II.

\begin{theorem}
Suppose that $F_0, F$ are continuous and
are quadratic forms of $x.$
If parameter-dependent S-procedure is lossless
then there exists a continuous function $h(.)$
that satisfies \eq{F0>hF}.
The set of all continuous  $h(.)$ satisfying
\eq{F0>hF} is open in $\cC(\pd, \bR^l).$
\end{theorem}

\sect{Generalized Yakubovich Lemma}

Let us begin from the formulation of the classical result.
\begin{theorem}[Yakubovich 1962, 1964]
\label{YL}
Let $A \in \bC^{n_x\times n_x},$
$B \in \bC^{n_x\times n_u},$ $G \in \bHM_{n_x + n_u}.$
The following statements are equivalent:

{\rm 1.} The inequality
\be
\vect{x}{u}^* G \vect{x}{u} > 0
\ee{xu*Gxu>0}
is fulfilled for all
$x \in \bC^{n_x\times 1}, \; u \in \bC^{n_u\times 1},
\; |x| + |u| \not= 0, \;
\omega \in [-\infty, +\infty]$ such
that
$\bi \omega x = A x + B u.$

{\rm 2.} There exists  $H \in \bHM_{n_x}$
that satisfies the LMI
\be
G > \matr{H A + A^* H}{H B}{B^* H}{0}.
\ee{G>HA+A*H}

\end{theorem}
Hereafter
$\bi$ is the imaginary unit.

This theorem was first proven by
V.A.Yakubovich in \cite{Yakubovich.1962}
for the single-input system (i.e. for $n_u = 1$) and in
\cite{Gantmaher.Yakubovich.1964} for the multi-input system.
Our formulation follows  \cite{Yakubovich.SMJ.1973}.
Now the statement of the Theorem~\ref{YL}
is known as a part of
Kalman-Yakubovich-Popov lemma
\cite{Yakubovich.1962,
Kalman.1963, Gantmaher.Yakubovich.1964, Popov.1964}.
In our opinion, it is more convenient for reference purposes,
and is more correct from the historical point of view,
to refer to Theorem~\ref{YL}
as Yakubovich Lemma. We call this statement a lemma to
point out its
connection with Kalman-Yakubovich-Popov Lemma. Additional
information about history of Kalman-Yakubovich-Popov Lemma can be
found in \cite{Gusev.Lihtarnikov.2006}.


The theorem below is proven in \cite{Iwasaki.Hara.AC.2005}.
It is a generalization of Yakubovich Lemma.
The paper \cite{Iwasaki.Hara.AC.2005} is based on generalizations of
Yakubovich Lemma that were proposed in
\cite{Churilov.1981,Iwasaki.Meinsma.Fu.2000}.

Let $\Theta = \Theta^* = \matr{\theta_{11}}{\theta_{12}}
{\theta_{21}}{\theta_{22}},$
$\det \Theta < 0.$
Define the curve
$\Gamma_{\Theta} = \{ \lambda \in \bar{\bC} \; | \;
(\lambda, 1) \Theta (\lambda, 1)^* = 0 \},$
and the domain
$\Omega_{\Theta} = \{ \lambda \in \bar{\bC} \; | \;
(\lambda, 1) \Theta (\lambda, 1)^* \ge 0 \}.$
Hereafter $\bar{\bC}$ is the closed complex plane
$\bar{\bC} = \bC \cup \{\infty\}.$
Choosing different matrices $\Theta$ it is possible to define any
circle or straight line on $\bar{\bC}.$
To illustrate, let us present some examples:
$\Theta = \matr{0}{1}{1}{0},
\Gamma_{\Theta} = \bi \bR \cup \infpoint,
\Omega_{\Theta} = \{ \Re \lambda \ge 0 \}\cup \infpoint;$
$\Theta = \matr{0}{-\bi}{\bi}{0},
\Gamma_{\Theta} = \bR \cup \infpoint,
\Omega_{\Theta} = \{ \Im \lambda \ge 0 \} \cup \infpoint;$
$\Theta = \matr{-1}{0}{0}{r^2}, \Gamma_{\Theta} = \{ | \lambda| = r\},
\Omega_{\Theta} = \{ |\lambda| \le r \}.$

Let $M, N \in \bC^{n\times n_z},$ $S \in \bHM_{n_z},$ $H \in \bHM_n.$
Define the generalized Lyapunov operator
$\Lambda_{M,N,\Theta}(S) = (M, N)(\Theta\otimes S)(M, N)^*.$
\sloppy
The adjoint operator takes the form
$\Lambda'_{M,N,\Theta}(H) = $
$(M^*, N^*)(\Theta^\top\!\otimes\! H)(M^*, N^*)^*.$
Hereafter $M_1\!\otimes\! M_2$ is a Kronecker product
of matrices $M_1$ and $M_2.$

\begin{theorem}[Iwasaki and Hara, 2005]
Let $M, N \in \bC^{n\times n_z},$ $G \in \bHM_{n_z}.$
Suppose  the intersection $\Gamma_{\Theta_1} \cap \Omega_{\Theta_2}$
includes more than one point and
$n_z > n;$
then the following statements are equivalent:

{\rm 1.} The inequality
\be
z^* G z > 0
\ee{z*Gz>0}
is fulfilled for all
$z \in \bC^{n_z\times 1}, z\not= 0,$
$\lambda \in \Gamma_{\Theta_1} \cap \Omega_{\Theta_2}$
such that
$(\lambda N - M) z = 0.$

{\rm 2.} There exist
$H_1, H_2 \in \bHM_n,$ $H_2 > 0$ such that
\be
G > \Lambda_{M,N,\Theta_1}'(H_1) + \Lambda_{M,N,\Theta_2}'(H_2).
\ee{G>Lambda12}

\end{theorem}


\sect{Yakubovich Lemma for matrix frequency-domain inequality}

Let $M, N, \Theta_i, i = 1,2$
be defined as in previous section.
Consider matrices $G_{ij} \in \bC^{n_z\times n_z},$
$G_{ij} = G^*_{ji}, \; i,j = 1, \ldots, m.$

\begin{theorem}
Suppose  the intersection $\Gamma_{\Theta_1} \cap \Omega_{\Theta_2}$
includes more than one point,
$n_z=n+1,$  and $\rank(\lambda N - M) = n$ for all
$\lambda \in \Gamma_{\Theta_1} \cap \Omega_{\Theta_2};$
then the following statements are equivalent:

{\rm 1.} The matrix inequality
\be \left(
\begin{array}{ccc}
z^* G_{11} z & \ldots & z^* G_{1m} z\\
\vdots & \ddots & \vdots \\
z^* G_{m1} z & \ldots & z^* G_{mm} z
\end{array}
\right) > 0
\ee{(Gii)>0}
is fulfilled for all
$z \in \bC^{n_z\times 1}, z \not= 0,$ $\lambda \in \Gamma_1 \cap \Omega_2$ such that
$(\lambda N - M) z = 0.$

{\rm 2.} There exist
$H_1, H_2 \in \bHM_{m n},$ $H_2 > 0$ such that
\be
\left(
\begin{array}{ccc}
G_{11} & \ldots &  G_{1m} \\
\vdots & \ddots & \vdots \\
 G_{m1}  & \ldots &  G_{mm}
\end{array}
\right) > \Lambda_{I_m\otimes M,I_m\otimes N,\Theta_1}'(H_1) +
\Lambda_{I_m\otimes M,I_m\otimes N,\Theta_2}'(H_2).
\ee{(Gii)>Lambda12}
\end{theorem}


\sect{Parameter-dependent Yakubovich Lemma}

Let $\pd$ be a topological space.
Consider maps
$G_{ij}(.) \in \cC(\pd, \bC^{n_z\times n_z}), \; $
$(G_{ij}(\pv) = G^*_{ji}(\pv) \;
\forall \pv \in \pd, \; i,j = 1, \ldots, m),$
$M(.), N(.) \in \cC(\pd, \bC^{n\times n_z}),$
$\Theta_i(.) \in \cC(\pd, \bHM_2), \;$
$(\det \Theta_i(\pv) < 0 \; \forall \pv \in \pd, \;
i = 1, 2).$
Denote $G(\pv) = \left(
\begin{array}{ccc}
G_{11}(\pv) & \ldots &  G_{1m}(\pv) \\
\vdots & \ddots & \vdots \\
 G_{m1}(\pv)  & \ldots &  G_{mm}(\pv)
\end{array}
\right).$

\begin{theorem}
Suppose $\pd$ is a compact space, the intersection
$\Gamma_{\Theta_1(\pv)} \cap \Omega_{\Theta_2(\pv)}$
includes more than one point for all $\pv \in \pd,$
and either $m=1, \; n_z > n$ or $m > 1,$
$n_z=n+1,$  $\rank(\lambda N(\pv) - M(\pv)) = n$ for all
$\lambda \in \Gamma_{\Theta_1(\pv)} \cap \Omega_{\Theta_2(\pv)}$
and all $\pv \in \pd;$
then the following statements are equivalent:

{\rm 1.} The matrix inequality
\be
\left(
\begin{array}{ccc}
z^* G_{11}(\pv) z & \ldots & z^* G_{1m}(\pv) z\\
\vdots & \ddots & \vdots \\
z^* G_{m1}(\pv) z & \ldots & z^* G_{mm}(\pv) z
\end{array}
\right) > 0
\ee{(Gii(p))>0}
is fulfilled for all
$\pv \in \pd,$ $z \in \bC^{n_z\times 1}, z\not= 0,$
$\lambda \in
\Gamma_{\Theta_1(\pv)} \cap \Omega_{\Theta_2(\pv)}$
such that
$(\lambda N(\pv) - M(\pv)) z = 0.$

{\rm 2.} There exist
$H_1(.), H_2(.) \in \cC(\pd, \bHM_{m n})$
such that for all $ \pv \in \pd$
\be
H_2(\pv) > 0,
G(\pv) > \Lambda_{I_m\otimes M(\pv),I_m\otimes N(\pv),\Theta_1(\pv)}'(H_1(\pv)) +
\Lambda_{I_m\otimes M(\pv),I_m\otimes N(\pv),\Theta_2(\pv)}'(H_2(\pv)).
\ee{(G(p))>Lambda12(p)}

The set of continuous pairs $(H_1(.),\! H_2(.))$ satisfying \eq{(G(p))>Lambda12(p)}
is open in $\cC(\pd, \bHM_{m n}\times \bHM_{m n}).$

\end{theorem}
\sect{Parameter-independent solutions of parameter-dependent LMI}

Consider the parameter-dependent LMI \eq{L>0}.
In contrast to section~2 now we are looking for
a constant vector $h$ that satisfies \eq{L>0}
for all $\pv \in \pd.$
In section~2 it was shown that when
$\cD \subset \bR^{k}$
the search of parameter-dependent solution of an LMI
can be reduced to the search of constant one
for a set of parameter-dependent LMIs  with growing
dimension of solution.
Now we would like to show that in some cases
the search
of the parameter-independent solution of
a parameter-dependent LMI can be reduced to
the solution of a set of parameter-independent LMIs
with growing dimension.

Suppose that $L$ is polynomial of $p$ and the set
$\pd$ is given by
\be
\pd = \{ \pv \in \bR^{1\times k} |
a_1 \le \pv_1 \le b_1, \;
a_i(p^{i-1}) \le \pv_i \le b_i(p^{i-1}),
\; i =  2, \ldots, k \},
\ee{D=}
where
$\pv = (p_1, \ldots, p_k),$
$\pv^i = (p_1, \ldots, p_i) \in \bR^{1\times i},$
$a_1 < b_1,$ $a_i, b_i, \; i =  2, \ldots, k,$
are polynomials and
$a_i(p^{i-1}) < b_i(p^{i-1})$ for all $p \in \pd.$

Let us describe the procedure for construction
of a parameter-independent LMI,
each solution of which defines the vector $h$ satisfying
\eq{L>0} for all $\pv \in \pd.$

Define $d = \lfloor\frac{\deg_{p_k}(L) + 1}{2}\rfloor,$
$\zeta(p_k) = (p_k^d, \ldots,  p_k, 1) \in \bR^{1\times(d+1)};$
then $L(p, h)$ can be represented as
$$L(p, h) = \left(
\begin{array}{ccc}
\zeta G_{11}(p^{k-1}, h) \zeta^* & \ldots & \zeta G_{1m}(p^{k-1}, h) \zeta^*\\
\vdots & \ddots & \vdots \\
\zeta G_{m1}(p^{k-1}, h) \zeta^* & \ldots & \zeta G_{mm}(p^{k-1}, h) \zeta^*
\end{array}
\right),$$
where the matrices
$G_{i j}(p^{k-1}, h), \; i, j = 1, \dots, m,$
are polynomials of $p^{k-1}$ and
are affine with respect to $h.$
For all $i, j = 1, \dots, m,$
$G_{i j}(p^{k-1}, h) = G_{j\,i}(p^{k-1}, h)^*.$

Define matrices
$M = (I_d, 0), \; N = (0, I_d) \in \bR^{d\times (d+1)},$
$$\Theta_1 = \matr{0}{-\bi}{\bi}{0},
\Theta_2(p^{k-1}) =
\matr{-1}{\frac{1}{2}(a_{k}(p^{k-1}) + b_{k}(p^{k-1}))}
{\frac{1}{2}(a_{k}(p^{k-1}) + b_{k}(p^{k-1}))}
{- a_{k}(p^{k-1}) b_{k}(p^{k-1})}.$$
Then
$$
(p_k N - M) z = 0 \Rightarrow \exists c \in \bC :
z = c (p_k^d,  \ldots, p_k, 1)^\top,
$$
$$\Gamma_{\Theta_1}\cap \Omega_{\Theta_2(p^{k-1})} =
[a_{k}(p^{k-1}),b_{k}(p^{k-1})].$$

Thus \eq{L>0} is fulfilled for all $p \in \cD$
iff
\be
\left( \begin{array}{ccc}
z^* G_{11}(p^{k-1}, h) z & \ldots & z^* G_{1m}(p^{k-1}, h) z\\
\vdots & \ddots & \vdots \\
z^* G_{m1}(p^{k-1}, h) z & \ldots & z^* G_{mm}(p^{k-1}, h) z
\end{array}
\right) > 0
\ee{zLz>0} 
for all  $p^{k-1} \in \cD^{k-1} , z \in \bC^{(d+1)\times 1}, z \not= 0,
p_k \in \Gamma_{\Theta_1}\cap \Omega_{\Theta_2(p^{k-1})},$
satisfying $(p_k N - M) z = 0.$
Here
$\cD^{k-1} =
 \{ q \in \bR^{1\times (k-1)} |
a_1 \le q_1 \le b_1, \;
a_i(q^{i-1}) \le q_i \le b_i(q^{i-1}),
\; i =  2, \ldots, k-1 \}.$

Let us introduce matrix polynomials
$H_1(p^{k-1}), H_2(p^{k-1}) \in \cC(\cD^{k-1}, \bHM_{d m}).$
Consider inequalities
\be
G(p^{k-1}, h) >
\Lambda_{I_m\otimes M,I_m\otimes N,
\Theta_1}'(H_1(p^{k-1})) +
\Lambda_{I_m\otimes M,I_m\otimes N,
\Theta_2(q)}'(H_2(p^{k-1})), \;
H_2(p^{k-1}) > 0.
\ee{Q>Lambda}
where
$$ G(p^{k-1}, h) =
\left( \begin{array}{ccc}
G_{11}(p^{k-1}, h)  & \ldots &  G_{1m}(p^{k-1}, h) \\
\vdots & \ddots & \vdots \\
 G_{m1}(p^{k-1}, h)  & \ldots &  G_{mm}(p^{k-1}, h)
\end{array}\right) \in \bHM_{m(d+1)}.
$$
From Theorem~6 it follows that if the vector $h$
and the pair of polynomials $H_1(p^{k-1}), H_2(p^{k-1})$
satisfies \eq{Q>Lambda} for all $p^{k-1} \in \cD^{k-1},$
then $h$ satisfies \eq{L>0} for all $p \in \cD.$

Let $h^{(1)}_{aux}$
be the vector of coefficients of polynomials
$H_1(p^{k-1}), H_2(p^{k-1});$
then inequalities \eq{Q>Lambda}
can be rewritten as follows:
\be
R_{1,i}(p^{k-1}, h, h^{(1)}_{aux}) > 0, \; i = 1,2,
\ee{R_1>0}
where
$ R_{1, 1}(p^{k-1}, h, h^{(1)}_{aux}) =
G(p^{k-1}, h) -
\Lambda_{I_m\otimes M,I_m\otimes N,\Theta_1}'(H_1(p^{k-1})) -
\Lambda_{I_m\otimes M,I_m\otimes N,\Theta_2(q)}'(H_2(p^{k-1})),$
$
R_{1, 2}(p^{k-1}, h, h^{(1)}_{aux}) = H_2(p^{k-1}).
$
The maps $R_{1, i}, \; i = 1,2,$
are polynomials of $p^{k-1}$ and are affine with respect
to the joint vector $(h, h^{(1)}_{aux}).$

Applying the same procedure to each inequality in
\eq{R_1>0}, we obtain
the system of four inequalities
\be
R_{2,i}(p^{k-2}, h, h^{(2)}_{aux}) > 0, \; i = 1,2,3,4,
\ee{R_2>0}
where $h^{(2)}_{aux}$ is the vector of parameters
of all polynomials, that were introduced
on two steps of procedure,
the maps $R_{2, i}, \; i = 1,2,3,4,$
are polynomials of $p^{k-2}$ and are affine with respect
to the joint vector $(h, h^{(2)}_{aux}).$

Repeating the procedure $k$-times,
we  obtain the system of
parameter-independent LMIs that has the form
\be
R_{k,i}(h, h_{aux}) > 0, \; i = 1, \dots , 2^k,
\ee{R_k>0}
where $h_{aux}$
is the vector of parameters
of all introduced polynomials,
the maps $R_{k, i}, \; i = 1, \dots , 2^k,$
are affine with respect
to the joint vector $(h, h_{aux}).$

For convenience we rewrite \eq{R_k>0}
as a single LMI
\be
R(h, h_{aux}) > 0,
\ee{R>0}
where
$R$ is block-diagonal matrix with blocks
$R_{k, i}, \; i = 1, \dots , 2^k,$
$h_{aux} \in \bR^{l_{R}}.$

Let  $\cR$ be
the set of all affine maps $R$ that can be constructed
from the map $L,$
using described recursive procedure.

\begin{theorem}
Suppose that $L$ is polynomial of $\pv,$
$\pd$ is
given by \eq{D=};
then
$$
\{ h \in \bR^l \; | \; L(p, h) > 0 \;
\forall p \in \cD \}
= \bigcup_{R \in \cR}
\{ h \in \bR^l \; | \;
\exists h_{aux} \in \bR^{l_R} : \; R(h, h_{aux}) > 0 \}.
$$
\end{theorem}


\sect{Example.
Minimization of a polynomial in a bounded domain}

Let $g(p)$ be the polynomial of $p \in \bR^{1\times k}.$
Consider the problem of evaluation of the  polynomial $g$
minimum in a domain $\cD$ given by \eq{D=}.

Define $L(p, h) = g(p) - h, \; h \in \bR.$
It is clear that any $h$ satisfying \eq{L>0}
for all $p \in \cD$ is a lower bound for $g$ in $\cD.$
Moreover from Theorem~7 it follows that
\be
\sup_{R \in \cR}\sup 
\{ h \; | \; \exists h_{aux} \in \bR^{l_{R}} : \; R(h, h_{aux}) > 0 \}
=
\min\{ g(p) \; | \; p \in \cD \}.
\ee{sup=min}
Taking  into account \eq{sup=min}, we can say that
minimization of the polynomial $g$ in a domain $\cD$
is reduced to a standard LMI optimization problem.

To illustrate the method we apply this approach
to a  polynomial of two variables $g(x, y).$
Let $\cD = \{ (x,y) \; |
\; a \le x \le b, c(x) \le y \le d(x) \},$
where $a < b,$ $c(x), d(x)$ are polynomials,
$c(x) < d(x)$ for all $x \in [a, b].$
Consider the parameter-dependent LMI in $h \in \bR$

\be g(x, y) - h > 0,  \forall x, y \in \cD
\ee{g-h>0}

Let $\deg_y g$ be the degree of the polynomial $g$
with respect to $y.$
Denote
$n_y = \lfloor\frac{\deg_{y}(g) + 1}{2}\rfloor,$
 $Y = (y^{n_y}, \ldots, y, 1).$
Then
$g(x, y) - h = Y G(x, h) Y^*,$
where \newline
$
G(x, h) = \left(\begin{array}{ccc}
G_{1, 1}(x) & \ldots & G_{1, n_y + 1}(x)\\
\vdots & \ddots & \vdots \\
G_{n_y + 1, 1}(x) & \ldots & G_{n_y + 1, n_y + 1}(x) - h
\end{array}
\right) \in \bHM_{n_y + 1},
$
$G_{ij}(x), i, j = 1, \ldots, n_y + 1,$ are polynomials of $x.$

Define matrices
$M_y = (I_{n_y}, 0), \; N_y = (0, I_{n_y})
\in \bR^{n_y\times (n_y+1)}.$
Then
\be
(y N_y - M_y) Z = 0
\ee{Zy}
implies $Z = c Y^\top,$ $c \in \bC.$

Define matrices
$$\Theta_1 = \matr{0}{-\bi}{\bi}{0},
\Theta_2(x) =
\matr{-1}{\frac{1}{2}(c(x) + d(x))}
{\frac{1}{2}(c(x) + d(x))}
{- c(x) d(x)}.$$
Then
$\Gamma_{\Theta_1}\cap \Omega_{\Theta_2(x)} =
[c(x), d(x)].$
Thus \eq{g-h>0}
is equivalent to
$Z^* G(y) Z > 0$ $\forall x \in [a, b],$
$y \in \Gamma_{\Theta_1}\cap \Omega_{\Theta_2(x)}$
and $Z \not= 0$
satisfying \eq{Zy}.

By Theorem~6 it follows that
the latter condition is fulfilled iff
there exist polynomial matrices
$H_{G1}(x),$ $H_{G2}(x) \in \cC(\bR, \bHM_{n_y})$
satisfying the LMIs
\be
H_{G2}(x) > 0, \; G(x, h) >
\Lambda'_{M_y, N_y, \Theta_1}(H_{G1}(x)) +
\Lambda'_{M_y, N_y, \Theta_2(x)}(H_{G2}(x))
\ee{G>L}

Let $n$ be the maximal degree of $x$ in \eq{G>L}.
Put $n_x = \lfloor\frac{n + 1}{2}\rfloor,$
$X = (x^{n_x}, \ldots, x, 1).$
Let us denote the
vectors of the coefficients of the polynomials
$H_{G1}(x), H_{G2}(x)$ by
$\hat{H}_{G1}, \hat{H}_{G2}$ respectively and
define matrices
$E(\hat{H}_{G2}, x) = H_{G2}(x),$
$F(\hat{H}_{G1}, \hat{H}_{G2}, x, h) = G(x, h) -
\Lambda'_{M_y, N_y, \Theta_1}(H_{G1}(x)) +
\Lambda'_{M_y, N_y, \Theta_2(x)}(H_{G2}(x)).$
Then LMIs \eq{G>L} can be rewritten as follows
\be
E(\hat{H}_{G2}, x) =
\left( \begin{array}{ccc}
X E_{1, 1}(\hat{H}_{G2}) X^* & \ldots
&X E_{1, n_y}(\hat{H}_{G2}) X^*\\
\vdots & \ddots & \vdots\\
X E_{n_y, 1}(\hat{H}_{G2}) X^* & \ldots
&X E_{n_y, n_y}(\hat{H}_{G2}) X^*
 \end{array} \right) > 0,
\ee{XEX*>0}
\be
F(\hat{H}_{G1}, \hat{H}_{G2}, x, h) =
\left( \begin{array}{ccc}
X F_{1, 1}(\hat{H}_{G1}, \hat{H}_{G2}) X^* & \ldots
&X F_{1, n_y+1}(\hat{H}_{G1}, \hat{H}_{G2}) X^*\\
\vdots & \ddots & \vdots\\
X F_{n_y+1, 1}(\hat{H}_{G1}, \hat{H}_{G2}) X^* & \ldots
&X F_{n_y+1, n_y+1}(\hat{H}_{G1}, \hat{H}_{G2}, h) X^*
 \end{array} \right) > 0.
\ee{XFX*>0}

Define matrices
$M_x = (I_{n_x}, 0), \; N_x = (0, I_{n_x})
\in \bR^{n_x\times (n_x+1)}.$
Then
\be
(x N_x - M_x) Z = 0
\ee{Zx}
implies $Z = c X^\top,$ $c \in \bC.$
Let $\Theta_3 =
\matr{-1}{\frac{1}{2}(a + b)}
{\frac{1}{2}(a + b)}
{- a b}.$

According to Theorem~5,
\eq{XEX*>0} and \eq{XFX*>0} are fulfilled
iff there exist $H_{E1}, H_{E2} \in \bHM_{n_x n_y},$
$H_{F1}, H_{F2} \in \bHM_{n_x (n_y + 1)}$
that satisfy the LMIs
\be
\begin{array}{c}
H_{E2} > 0, \; H_{F2} > 0, \; 
\left( \begin{array}{ccc}
E_{1, 1}(\hat{H}_{G2})  & \ldots &
E_{1, n_y}(\hat{H}_{G2}) \\
\vdots & \ddots & \vdots \\
E_{n_y, 1}(\hat{H}_{G2})  & \ldots &
E_{n_y, n_y}(\hat{H}_{G2})
\end{array}
\right) >
\begin{array}{l}
\\
\Lambda_{I_{n_y}\otimes M_x,I_{n_y}\otimes N_x,\Theta_1}'
(H_{E1}) + \\
\Lambda_{I_{n_y}\otimes M_x,I_{n_y}\otimes N_x,\Theta_3}'
(H_{E2}),
\end{array}\\ \\
\left( \begin{array}{ccc}
F_{1 1}(\hat{H}_{G1}, \hat{H}_{G2})  & \ldots &
F_{1 n_y+1}(\hat{H}_{G1}, \hat{H}_{G2}) \\
\vdots & \ddots & \vdots \\
F_{n_y+1 1}(\hat{H}_{G1}, \hat{H}_{G2})  & \ldots &
F_{n_y+1 n_y+1}(\hat{H}_{G1}, \hat{H}_{G2}, h)
\end{array}
\right) >
\begin{array}{l}
\\
\Lambda_{I_{n_y+1}\otimes M_x,I_{n_y+1}\otimes N_x,\Theta_1}'(H_{F1}) + \\
\Lambda_{I_{n_y+1}\otimes M_x,I_{n_y+1}\otimes N_x,\Theta_3}'(H_{F2}).
\end{array}
\end{array}
\ee{finLMIs}

In this way \eq{g-h>0} is fulfilled iff
there exists $n_x$ such that LMIs \eq{finLMIs} are feasible.
Evidently, system \eq{finLMIs}  can be written as a single
LMI \eq{R>0}.
Maximizing $h$ over solutions of \eq{finLMIs} and
increasing $n_x$ we can estimate minimum of $g$ in $\cD$
with any desired accuracy.

The proposed approach can be applied to
verification of positivity of polynomials.
As an illustration let us consider a numerical example.

It is known that
$g(x, y) = x^4 y^2 + x^2 y^4 -3 x^2 y^2 + 1 +
\varepsilon  > 0$
for all $ \varepsilon > 0$ and all $x, y.$
Besides for any $\varepsilon$ the polynomial $g(x, y)$
cannot be represented as
sum of squares (SOS) of
polynomials~\cite{Parillo.thesis.2000}.
It means that the well known SOS representation technique
cannot be used directly to verify the positivity of $g.$

It is easy to see that $g(x, y) > 0,$ when $|x| >2$ or $|y| > 2.$
So, to verify the positivity of $g$ it is sufficient to
find positive lower bound of $g(x,y)$ in
$\cD = \{ (x, y) \; | \; |x| \le 2, |y| \le 2\}.$
Calculations using Matlab$^\circledR$ LMI toolbox showed that
for $n_x = 2$ the system of inequalities \eq{finLMIs}
has a solution with $h > 0$ if
$\varepsilon = 10^{-11}.$
This proves that
$x^4 y^2 + x^2 y^4 -3 x^2 y^2 + 1 + 10^{-11} > 0$
for all $x, y.$


\sect{Construction of parameter-dependent
Lyapunov function}

Consider the parameter-dependent nonlinear system
\be
\dot{x} = A(p) x + B(p) u, \; x(0) = x_0,
\ee{dx=Ax+Bu}
\be
u = \varphi(p, t, x),
\ee{u=f}
where $t \ge 0,$ $x, x_0 \in \bR^{n\times 1},
u \in \bR^{m\times 1},$ $p$ is a parameter,
$p \in \cD \subset \bR^{1\times k},$
$A(.) : \cD \to \bR^{n\times n},$
$B(.) : \cD \to \bR^{n\times m},$
$\varphi : \cD\times [0, +\infty)\times \bR^{n\times 1}
\to \bR^{m\times 1}.$

The nonlinearity $\varphi$ satisfies the quadratic constraint
\be
(x^*, u^*)G(p)(x^*, u^*)^* \ge 0
\ee{xuGxu>0}
which is fulfilled for all $x \in \bR^{n\times 1},
t \ge 0, p \in \cD, u = \varphi(p, t, x).$
Here $G(p) \in \bHM_{m+n}.$

Suppose that $A, B,$ and $G$ are polynomials of $p,$ and
$\cD$ is defined by \eq{D=}.
Consider the parameter-dependent Lyapunov function candidate
$V(p, x) = x^* H(p) x,$ where $H(.) \in \cC(\cD,\bHM_n).$
We are looking for $H(.)$ that satisfies
\be
H(p) > 0 \; \; \forall p \in \cD,
\ee{H>0}
and
\be
\frac{d}{dt} V(p, x(t)) < 0 \; \; \forall t \ge 0,
\ee{dV<0}
for all solutions of \eq{dx=Ax+Bu} satisfying \eq{xuGxu>0}.
If fulfilled these conditions guarantee the
asymptotic stability of the closed-loop system
\eq{dx=Ax+Bu}, \eq{u=f} for any $\varphi$
such that \eq{xuGxu>0} holds.

Define matrices
$M(p) = (A(p), B(p)), N = (I_n, 0) \in \bR^{n\times (n+m)},$
$\Theta = \matr{0}{1}{1}{0}.$
Then $$\Lambda'_{M(p), N, \Theta}(H(p)) =
\matr{H(p) A(p) + A^*(p) H(p)}{H(p) B(p)}
{B^*(p) H(p)}{0}
$$
and condition \eq{dV<0} takes the form:
\be
\forall p \in \cD,
\; \forall |x| +|u| \not= 0
\mbox{~\eq{xuGxu>0}~implies~}
(x^*, u^*)\Lambda'_{M(p), N, \Theta}(H(p))(x^*, u^*)^* < 0
\ee{I}
Using parameter-dependent S-procedure and Theorem~1
we can see that
\eq{I} is fulfilled
iff there are a polynomial matrix $H(p)$
and a polynomial $\eta(p)$
that satisfy
\be
\Lambda'_{M(p), N, \Theta}(H(p)) + \eta(p) G(p) >
0, \; \; \eta(p) > 0
\; \forall p \in \cD.
\ee{()+G>0}

Let $h \in \bR^{l}$ be the vector of coefficients of polynomials
$H(p), \eta(p).$
Then \eq{()+G>0} can be written as a single
parameter-dependent LMI \eq{L>0}.
Using the procedure described in Section~7  we can define
the set of affine maps $\cR.$
From Theorem~7 it follows that
$H(p)$ and $\eta(p)$
satisfying \eq{H>0} and \eq{()+G>0}
exists iff there is  an affine map $R \in \cR$
such that  LMI \eq{R>0} is feasible.
Any solution $h$ of  obtained in this way LMI \eq{R>0}
defines the polynomial
matrix $H(p)$ that satisfies \eq{H>0} and \eq{dV<0}.\\


\end{document}